\documentclass[11 pt]{article}    
\usepackage[margin = 1in]{geometry}
\usepackage{amsfonts, amssymb}					 
\usepackage[colorlinks=true,citecolor=teal,linkcolor=blue]{hyperref}
\usepackage{tikz}
\usepackage{tikz-cd}
\usepackage{quiver} 
\usepackage{xcolor}
\usepackage{graphicx} 
\usepackage{amsmath}
\usepackage{mathrsfs}
\usepackage{mathabx}
\usepackage{comment}
\usepackage{setspace} 
\usepackage{enumerate}
\usepackage{amsthm}

\usepackage{url}

\newcommand{\eqn}[0]{\begin{array}{rcl}}
\newcommand{\eqnend}[0]{\end{array} }  	
\newcommand\norm[1]{\left\lVert#1\right\rVert}
\newcommand\restr[2]{{
  \left.\kern-\nulldelimiterspace 
  #1
  \vphantom{\big|} 
  \right|_{#2} 
  }}
  \newcommand{\inn}[1]{\langle#1\rangle}
\def\P{\mathbb{P}}

\def\N{\mathbb{N}}
\def\R{\mathbb{R}}

\numberwithin{equation}{section}
\begin{document}	

\newtheorem{theorem}{Theorem}[section]
\newtheorem{prop}[theorem]{Proposition}

\theoremstyle{definition}
\newtheorem{definition}[theorem]{Definition}

\newtheorem{assumption}[theorem]{Assumption}

\theoremstyle{definition}
\newtheorem{remark}[theorem]{Remark}

\theoremstyle{definition}
\newtheorem{problem}[theorem]{Problem}
\newtheorem{conjecture}[theorem]{Conjecture}

\theoremstyle{definition}
\newtheorem{example}[theorem]{Example}

\newtheorem{corollary}[theorem]{Corollary}

\newtheorem{lemma}[theorem]{Lemma}
\theoremstyle{definition}
\newtheorem{notation}[theorem]{Notation}
\title{Spectral Rigidity of Stochastic Operators with Heavy Tailed Perturbations}
\author{Youssef Djellouli\footnote{Syracuse University, Department of Mathematics. Email:
\href{mailto:mydjello@syr.edu}{mydjello@syr.edu}.}}
\date{}
\maketitle

\begin{abstract}
In this paper, we provide a new set of tools for establishing the number rigidity of point processes which are particularly well suited for the spectral processes of stochastic operators. These tools allow one to prove the number rigidity of a wide class of stochastic operators that had been previously inaccessible to the standard methods in the area, such as the variance of linear statistics, which require being able to precisely analyze correlations between individual points. We propose an alternative mechanism which uses a pathwise comparison between the growth rate of the unperturbed eigenvalues and the cumulative random noise which allows one to convert the probabilistic problem into a more tractable perturbation theory problem. This new approach allows us to prove the rigidity of previously inaccessible operators including those with heavy tailed perturbations or complex correlation structures. As a demonstration, we use these tools to prove the number rigidity of new classes of stochastic Schr\"odinger and Sturm-Liouville Operators which are resistant to analysis using the standard rigidity toolkit.

\end{abstract}

\section{Introduction}
Let $\Pi$ be a locally finite point process on $\R$. Informally, $\Pi$ is said to be number rigid if one can recover the number of points inside a bounded subset of $\R$ using information from the point configuration outside of the subset. Number rigid processes contrast significantly with processes like the Poisson point process which possess the property that the number of points in disjoint sets are independent of each other. 

Since being introduced in \cite{GhoshPeres,HolroydSoo}, number rigidity has been studied for a large class of point processes. For instance, rigidity and related properties have been studied for various determinantal and Pfaffian point processes (e.g \cite{GhoshKrishnapur,Buf1,Buf2,Buf4,BufetovAiry,Heisenberg}), the eigenvalue processes of random matrices (e.g. \cite{ChhaibiNajnudel,Ghosh,RieszGas,Sine-BetaDLR,RigidityMSAO} and many others), processes of random operators arising from statistical physics ensembles (e.g. \cite{Sine-BetaDLR,RieszGas, LamarreGhosalLiaoSpacialConditioning,LamarreGhosalLiaoSpectralRigidity,RigidityMSAO,ChhaibiNajnudel,Leble}), lattice models (e.g. \cite{PeresSly,GhoshKrishnapur}) and
the spectral process of random Schr\"odinger operators (e.g. \cite{LamarreGhosalLiaoSpacialConditioning,LamarreGhosalLiaoSpectralRigidity,RigidityMSAO,BufetovAiry}). The eigenvalue processes of random Schr\"odinger operators are of particular interest in this paper. Indeed, there is a natural motivation for studying these processes stemming from statistical mechanics wherein the eigenvalues of the operator can be interpreted as the energy levels present in a quantum mechanical system with some disorder or impurity. 

In this paper, which is a follow up to \cite{djellouli2025rigidityprojectedperturbedlattices}, we extend the main results of \cite{djellouli2025rigidityprojectedperturbedlattices} and further develop new convenient tools for proving number rigidity which are particularly well suited for the spectral processes of stochastic operators. Then we apply these tools to new classes of stochastic Schr\"odinger and Sturm-Liouville operators with heavy tailed noises that have previously been inaccessible when using the standard rigidity toolkit. 

\subsection{Motivation and Set-Up}

\begin{definition}
We say that a locally finite point process $\Pi$ on $\R$ is number rigid if for every bounded Borel set $B\subset\R$, there exists a measurable function $F$ such that $F(\Pi\cap B^c)=|\Pi\cap B|$ almost surely.
\end{definition}

\begin{definition}
For a point process $\Pi$, and $S\subset\Pi$ finite, $S\neq\varnothing$, let \[\Pi_S=\Pi-\sum_{x\in S}\delta_{x}\]
Then, we say $\Pi$ is deletion singular if for every $S$ the laws of $\Pi$ and $\Pi_S$ are singular.
\end{definition}

The links between deletion singularity and number rigidity have been well studied (see e.g. \cite{PeresSly,HolroydSoo,GhoshPeres,Palm}). It is well known that for most typical point processes, number rigidity implies deletion singularity although the converse is not true in general. \\

In \cite{GhoshPeres}, Ghosh and Peres introduce a general tool for proving the number rigidity and deletion singularity of point processes by using the variance of linear statistics:
\begin{prop}\label{prop: variance of linear statistics}
    If there exists a sequence of Borel measurable functions $f_n:\R\to\R$ such that $f_n\to1$ uniformly on compact sets and 
    \[\lim_{n\to\infty}\mathrm{Var}\left[\sum_{x\in\Pi}f_n(x)\right]=0\]
    Then $\Pi$ is number rigid and deletion singular.
\end{prop}
This tool has been crucial for a very successful research program to prove the rigidity of many processes originating from a wide swathe of probability and mathematical physics. Some examples include lattice models (see e.g. \cite{GhoshKrishnapur}), random Schr\"odinger Operators \cite{LamarreGhosalLiaoSpacialConditioning,LamarreGhosalLiaoSpectralRigidity,RigidityMSAO}, the zeroes process of analytic functions and the spectral processes of scaled matrix ensembles (e.g. \cite{GhoshPeres}), and many other processes (see e.g. \cite{BufetovAiry,Buf4,Ghosh,StealthyHyperuniformity} and many others). 

Despite its diverse applications, the variance of linear statistics is not universally applicable. Proposition \ref{prop: variance of linear statistics} can be readily applied to processes with i.i.d noise since we have 
\[\mathrm{Var}\left[\sum_{x\in\Pi}f_n(x)\right]=\sum_{x,y\in\Pi}\mathrm{Cov}[f_n(x),f_n(y)]=\sum_{x\in\Pi}\mathrm{Var}(f_n(x))\]
However, it is apparent that this quantity becomes very difficult to control and is often completely intractable if the process contains non-trivial correlations. In particular, when considering the spectral processes of stochastic operators, this amounts to needing precise control over the correlations of individual eigenvalues which is often very difficult to achieve for many processes of interest.

\begin{remark}\label{rmk: exponential semigroup}
A particularly effective approach for investigating the number rigidity of the eigenvalue process of the stochastic operators we are interested in emerged from using the trace of the operator's semigroup in concert with the variance of linear statistics (see e.g. \cite{LamarreGhosalLiaoSpectralRigidity,LamarreGhosalLiaoSpacialConditioning,RigidityMSAO}). In particular, the points in the eigenvalue process of a stochastic operator $H$ are precisely its (random) eigenvalues $(\lambda_n)_{n\in\N}$.
Then, as is typical in applications of Proposition \ref{prop: variance of linear statistics}, one would like to use exponential test functions of the form \[f_n(\lambda) = e^{-t_n\lambda}\] for some sequence $t_n\to 0$ and so the natural question is how one can compute $\mathrm{Var}[{\sum_k}e^{-t_n\lambda_k}]$. A common way to do this is to use the trace of the exponential semigroup of $H$ which is defined by\[\mathrm{Tr}[e^{-tH}] = \sum_{k=1}^\infty e^{-t\lambda_k}\]
Then, the variance of this trace can in principle be computed by employing Feynman-Kac formulas for the semigroup provided that $H$ is sufficiently regular. Notice that this allows one to bypass the issue of controlling the correlations between individual eigenvalues. 
\end{remark}

The approach outlined in Remark \ref{rmk: exponential semigroup} cannot be applied to the eigenvalue process of every stochastic operator of interest. Indeed, there are some very important processes, such as the Airy-$2$ process, for whom the variance of its exponential semigroup trace is bounded above $0$ but which are nonetheless known to be rigid (see \cite[Proposition 2.27]{LamarreGhosalLiaoSpectralRigidity}). In the case of the Airy-$2$ process number rigidity was shown by Bufetov in \cite{BufetovAiry} by using the specific determinantal structure of that process in the $\beta=2$ case. However this approach does not generalize to other operators or even other members of the $\beta$-Airy family which do not share the determinantal structure. The semigroup trace method also fails in the case that the eigenvalue perturbations are heavy tailed since the exponential trace moments may be infinite in this case. As such there is a very natural motivation to search for other statistics that sidestep both the heavy tailed issue of the exponential semigroup and more generally the difficulty in handling complex correlations that comes from using  Proposition \ref{prop: variance of linear statistics} in the hopes of gaining a more general and convenient technique for proving the rigidity of spectral processes. 

\subsection{New Results}

In the paper \cite{djellouli2025rigidityprojectedperturbedlattices}, the authors propose a family of point processes emblematic of those whose rigidity is not amenable to being proven using the variance of linear statistics called the ``Projected Perturbed Lattice". They then develop a new scheme for proving deletion singularity for this family as well as more general point processes that can be adapted to those with heavy tailed perturbations or noise with complicated correlations. In particular, the method in \cite{djellouli2025rigidityprojectedperturbedlattices} does not rely on the existence of any moments and instead capitalizes on identifying ``mismatches" between the process and a sort of ``preimage". Then by comparing the worst case mismatch with the growth of said preimage one is able to determine whether points had been deleted from the process. Notably, this approach also does not require any independence of the noise (see Theorem \ref{thm: original result} for more details).\\ 

In this paper, we are interested in point processes of the form 
\[\Pi = \sum_{n\in I}\delta_{r_n+g_n}\] 
where we let $I$ be some countable index set, $(r_n)_{n\in I}\subset\R$ is unbounded either from above or below or both, and $(g_n)_{n\in I}$ are random variables taking values in $\R$ such that $\Pi$ is almost surely locally finite.
We think of the set $(r_n)_{n\in I}$ as the ''starting locations" of the process that are then perturbed by $(g_n)_{n\in I}$.
\begin{remark}\label{rem: operator set-up}
The eigenvalue process of a stochastic operator is just a a special case of the process $\Pi$. In particular, let $H$ be some (deterministic) self-adjoint differential operator with compact resolvent. Denote its spectrum by $(r_n)_{n\in I}$. Then let $\hat{H}$ be some perturbation of $H$, say by adding a random potential, that is also self-adjoint with compact resolvent, with spectrum $(\lambda_n)_{n\in I}$. Then the points $\lambda_n$ are exactly of the form $r_n+g_n$ if we let $g_n:=\lambda_n-r_n$. Note that we allow for repeated eigenvalues in both $(r_n)_{n\in I}$ and $(\lambda_n)_{n\in I}$.
\end{remark}

The main result in \cite{djellouli2025rigidityprojectedperturbedlattices} as applied to the model $\Pi$ can be summarized as follows:

\begin{theorem}\cite[Theorem 2.2]{djellouli2025rigidityprojectedperturbedlattices}:\label{thm: original result}
Enumerate $(r_n)_{n\in I}$ as $-\infty<r_1\leq r_2\leq\cdots\leq r_i\leq\cdots<\infty$ and suppose that the sequence $(r_n)$ is unbounded from above and bounded from below with no accumulation points. Let $(r_{n_k})_{n_k\in I}$ be the subsequence of $(r_n)_{n\in I}$ that includes every unique $r_i$ and the largest index of any repeated $r_i$'s. Finally, let $M_n:=\displaystyle\max_{i\leq n}|g_i|$. Then $\Pi$ is deletion singular if
\begin{align}\label{eqn: gap condition}
    \lim_{k\to\infty}\frac{M_{n_k}}{r_{n_k}-r_{n_{k-1}}}=0 \quad \text{a.s}.
\end{align}
\end{theorem}

\begin{remark}\label{rmk: original theorem is stronger}
    In fact the result proven in \cite{djellouli2025rigidityprojectedperturbedlattices} is even stronger and gives us that the laws of $\Pi_S$ and $\Pi_T$ are singular whenever $|S|\neq|T|$ i.e. one is able to deduce the precise number of deleted points. The argument is also valid under similar assumptions in the cases that $(r_n)$ is bounded from above and unbounded from below, or unbounded in both directions.

\end{remark}

In this paper, we extend Theorem \ref{thm: original result} and prove that number rigidity holds in addition to deletion singularity under similar and more convenient conditions that are more accessible to the standard mathematical physics toolkit (e.g. Weyl Laws). We then use these new conditions to prove the number rigidity of several operators of physical interest which heretofore have been inaccessible to the standard rigidity tools. \\

We begin with a statement that gives us number rigidity from a new condition that is more convenient to use for the perturbed operator setting we are interested in. 

\begin{prop}\label{prop: strong sum condition result}
Enumerate $(r_n)$ as $-\infty<r_1\leq r_2\leq\cdots\leq r_i\leq\cdots<\infty$ and suppose that the sequence $(r_n)_{n \in I}$ is unbounded from above and bounded from below with no accumulation points. Suppose additionally that 
\begin{align}\label{eqn: strong sum conditions}
    \sum_{i=1}^n|g_i| = \mathrm{o}(r_n) \text{ almost surely} \quad \text{as} \quad n\to\infty
\end{align}
Then $\Pi$ is number rigid.
\end{prop}
\begin{proof}
    See Section \ref{sec: proof of strong result}
\end{proof}

The primary benefit of this approach to number rigidity comes from the novel pathwise mechanism, \eqref{eqn: strong sum conditions}, which contrasts greatly with the prior moment methods of Proposition \ref{prop: variance of linear statistics} and others (e.g. \cite{HolroydSoo,PeresSly}). Rather than considering the variance of the perturbed eigenvalues in isolation, Proposition \ref{prop: strong sum condition result} instead compares the cumulative growth of the perturbations with the growth of the original unperturbed eigenvalues. This new statistic allows us to completely sidestep the problem of heavy tailed perturbations that can render the approach outlined in Remark \ref{rmk: exponential semigroup} completely impossible and (as we demonstrate in Section \ref{sec: applications}) also allows for considerably simpler computations by giving us access to a wide range of powerful techniques from perturbation theory. To the best of our knowledge, this is the first general method for proving number rigidity for operators with heavy tailed perturbations.

Additionally, we have the following corollary to Proposition \ref{prop: strong sum condition result} which trades a slightly stronger condition for some additional practicality in the perturbed operator setting.

\begin{corollary}\label{prop: weak sum condition result}
Enumerate $(r_n)$ as $-\infty<r_1\leq r_2\leq\cdots\leq r_i\leq\cdots<\infty$ and suppose that the sequence $(r_n)$ is unbounded from above and bounded from below with no accumulation points. Suppose additionally that 
\begin{align}\label{eqn: weak sum conditions}
    \sup_{n\geq1}|g_n|<\infty \quad \text{and} \quad \lim_{n\to\infty}\frac{r_n}{n}=\infty \text{ almost surely}
\end{align}
Then $\Pi$ is number rigid.
\end{corollary}
\begin{proof}
Corollary \ref{prop: weak sum condition result} follows immediately from Proposition \ref{prop: strong sum condition result}. Indeed it is easy to see that the conditions \eqref{eqn: weak sum conditions} directly implies the weaker condition \eqref{eqn: strong sum conditions} from the bound 
\[ \sum_{i=1}^n|g_i|\leq n\sup_{i\leq n}|g_i|\]
\end{proof}

\begin{remark}\label{rmk: positivity remark}
Similarly to the discussion in Remark \ref{rmk: original theorem is stronger}, both Propositions \ref{prop: weak sum condition result} and \ref{prop: strong sum condition result} can be easily adapted to account for $(r_n)_{n\in I}$ being unbounded in the negative direction. We omit these details because they are unnecessary for the applications we are interested in. In particular, it is easy to see that proving the rigidity of the eigenvalue process of an operator $H$ is equivalent to proving the rigidity of $-H$, so there is no loss in only considering the case where the unperturbed eigenvalues $(r_n)_{n\in I}$ are unbounded in the positive direction as opposed to the negative.
\end{remark}

\begin{notation}
We say that $f(x)\lesssim g(x)$ if there exists a constant $c>0$ such that $f(x)\leq cg(x)$ for every $x$. We say $f(x)\asymp g(x)$ if $f(x)\lesssim g(x)$ and $g(x)\lesssim f(x)$.
\end{notation}

\section{Applications}\label{sec: applications}
As outlined in the introduction, our primary motivation is to provide a convenient toolkit to prove the number rigidity of the eigenvalue processes of a variety of stochastic operators of independent physical and mathematical interest that have resisted prior methods.

\subsection{Sturm-Liouville Operator}
Consider the Sturm-Liouville Operator given by
\[Hu = \frac{1}{w(x)}\left(-\frac{d}{dx}(p(x)u'(x))+q(x)u(x)\right)\text{ on the Hilbert Space } L^2((a,b),w(x)dx)\]
with the standard regularity and self-adjoint boundary conditions (e.g. Dirichlet) and $p,w>0$ a.e and $p^{-1},w,q\in L^1(a,b)$. Then $H$ has compact resolvent and we let $r_1\leq r_2\leq...<\infty$ denote the spectrum of $H$. 

We define the perturbed operator $\hat{H}(\omega)$ by 
\[\hat{H}(\omega)u = \frac{1}{w(x)}\left(-\frac{d}{dx}(p(x)u'(x))+(q(x)+\xi_\omega(x))u(x)\right) \] with $\xi_\omega(x)$ random, satisfying $\frac{\xi_\omega(x)}{w}\in L^\infty(a,b)$ almost surely (preserving the regularity of $\hat{H}$) and let $\lambda_1(\omega)\leq\lambda_2(\omega)\leq\cdots<\infty$ denote the spectrum of $\hat{H}(\omega)$

It is well known (see e.g. \cite{SturmLiouville} for a comprehensive treatment) that \[r_n\asymp n^2 \text{ as } n\to\infty\] and by an elementary min-max argument we have 
\[|g_n(\omega)|:=|\lambda_n(\omega)-r_n|\leq\norm{\frac{\xi}{w}}_\infty\]

Thus it follows immediately from Corollary \ref{prop: weak sum condition result} that the eigenvalue point process $\{\lambda_n(\omega)\}_{n\in\N}$ is number rigid. 

Note that we do not need the noise to be stationary or independent. Indeed, one can take $\frac{\xi_\omega(x)}{w} = -XU(x)$ where $U(x)\in L^\infty$ and $U(x)\geq0$, $U$ not identically $0$, where $X$ is a random variable has no positive exponential moments (for example $X$ can be Pareto distributed). In this case, the semigroup trace approach outline in Remark \ref{rmk: exponential semigroup} becomes impossible because the moments of $e^{-t\lambda_n}=e^{-t(r_n+g_n)}$ will be infinite. 
However, we are nonetheless able to prove that the process is number rigid by using the criterion \eqref{eqn: weak sum conditions}.

\subsection{Schr\"odinger Operator with Poisson Shot Noise}

Consider the Schr\"odinger Operator given by 
\[H = -\frac{d^2}{dx^2}+V(x) \text{ on } \R \text{ with } V(x)\asymp|x|^\alpha\]
The standard Weyl Law gives us 
\[r_n\asymp n^{\frac{2\alpha}{\alpha+2}} \quad \text{as } n\to\infty\]
It is well known that this holds for all sufficiently nice potentials $V$ (see \cite{WeylLaw} for more details).
In particular, we have
\[\lim_{n\to\infty}\frac{r_n}{n}=\infty \iff \alpha>2\]

Consider the Poisson Shot Noise given by 
\[W(x) = \sum_{i=1}^{N(\omega)}A_i\phi(x-X_i)\]
where $\{X_i\}$ is a homogeneous Poisson process, $N<\infty$ almost surely, $\phi$ is compactly supported with $\phi\in L^1\cap L^\infty$, and each $A_i\geq0$ is heavy tailed and almost surely finite. 

Then we immediately have that the operator 
\[\hat{H}(\omega) = H - W\]
is number rigid when $\alpha>2$ by Corollary \ref{prop: weak sum condition result}.

Like the previous Sturm-Liouville example, the semigroup trace can have infinite moments rendering the standard approach described in Remark \ref{rmk: exponential semigroup} impossible; nonetheless we are able to obtain number rigidity for this operator through Proposition \ref{prop: strong sum condition result}

\subsection{Schr\"odinger Operator with Stable Perturbations}
Here we illustrate one of the major strengths of our new approach by considering Schr\"odinger Operators with heavy tailed noises such as Brownian motions and stable process perturbations. In particular, we show how relatively simple perturbation theory calculations can be translated into number rigidity for operators which were inaccessible with the previous rigidity toolkit.\\

For $m>0$ consider the Schr\"odinger Operator on $L^2(D)$ where $D$ is either $\R$ or $\R_+$, with Dirichlet boundary conditions given by 
\[H = -\Delta+|x|^m\] and the perturbed Schr\"odinger Operator given by 
\[\hat{H} = H+g(x)\]
where $g(x)$ is a centered, symmetric $\alpha$-stable Levy Process for $\alpha\in (0,2]$ (recall that when $\alpha=2$ this process is simply a Brownian motion).

\begin{prop}
$\hat{H}$'s eigenvalue process number rigid for every $m>2+\frac{2}{\alpha}$. 
\end{prop}
\begin{proof}
Our proof employs a perturbation theory argument that allows us to bound the gaps between the eigenvalues of the original operator $H$ and those of the perturbed operator $\hat{H}$  combined with estimates on the asymptotic (pathwise) growth of the $\alpha$-stable Levy Processes.\\

Suppose $g(x)$ is a real valued and locally bounded function satisfying:
\begin{align}\label{eqn: g bound}
    |g(x)|\leq C(1+h(|x|))
\end{align}
where $h(x)$ is a real-valued, eventually non-decreasing function that satisfies
\begin{align}\label{eqn: h conditions}
    h(x) = \mathrm{o}(x^m) \text{ and } \frac{h(x)}{x^m} \text{ is eventually non-increasing}
\end{align}
As usual, let $(r_n)_{n\in\N}$ denote the eigenvalues of original operator $H$ and $(\lambda_n)_{n\in\N}$ denote the eigenvalues of the perturbed operator $\hat{H}$.\\
Thanks to \eqref{eqn: h conditions}, if we take $R>0$ to be sufficiently large and split the domain into $|x|>R$ and $|x|\leq R$ we attain the bound
\begin{align}\label{eqn: g growth bound}
    |g(x)|\leq C +C'h(R)+C'\frac{h(R)}{R^m}|x|^m
\end{align}
Let 
\[q(u) = \inn{H^{1/2}u,H^{1/2}u}\] denote the quadratic form associated with the unperturbed operator $H$ and 
\[v(u) = \int_D g(x)|u(x)|^2dx\]
denote the quadratic form associated with the perturbation by $g(x)$.

Now using \eqref{eqn: g growth bound} gives us
\begin{align*}
    |v(u)|\leq \int_{|x|\leq R}|g(x)||u(x)|^2dx+\int_{|x|>R} |g(x)||u(x)|^2 dx
\end{align*}
For the first integral, by the monotonicity of $h(x)$ and \eqref{eqn: g bound} we have 
\[\int_{|x|\leq R}|g(x)||u(x)|^2dx\leq C(1+h(R))
\norm{u}_2^2\]
Then for the second integral, \eqref{eqn: h conditions} give us
\[|g(x)|\leq C+C\frac{h(R)}{R^m}|x|^m\] thus we get
\[\int_{|x|>R} |g(x)||u(x)|^2 dx \leq C\norm{u}_2^2+C\frac{h(R)}{R^m}\int_{D} |x|^m|u(x)|^2 dx\]
However, this final integral clearly is bounded by the form $q$, i.e.
\[\int_{D} |x|^m|u(x)|^2 dx\leq q(u)\] and so we conclude

\begin{align}\label{eqn: perturbation form bound}
    |v(u)| 
\leq C(1+h(R))\norm{u}_2^2+C\frac{h(R)}{R^m}q(u)
\end{align}

Hence, by taking $R:=r_n^{1/m}$ to be sufficiently large we are able to conclude from \cite[Theorem 4.13]{Veselic} that 
\begin{align}\label{eqn: singular eigenvalue spread}
    |\lambda_n-r_n|\leq C+2C'h(r_n^{1/m})
\end{align}

Indeed, we see that $\hat{H}$ satisfies the necessary conditions for Theorem 4.13 owing to the fact that $\hat{H}$ has no essential spectrum, its eigenvalues are bounded from below, and the bound \eqref{eqn: perturbation form bound} ensures that condition (ii) of the theorem is satisfied. Notice that the conditions for Theorem 4.13 require that the term $C\frac{h(R)}{R^m}$ in the bound \eqref{eqn: perturbation form bound} be strictly less than $1$, but this can of course be achieved by taking $R$ to be sufficiently large.\\

With the perturbation bound \eqref{eqn: singular eigenvalue spread} in hand, recall that from the classical Weyl Law we have 
\[r_n\asymp n^{\frac{2m}{m+2}} \text{ as } n\to\infty\]

First consider the case where $g(x)=B(x)$ is Brownian motion (this corresponds to the $\alpha=2$ case).
Then by the law of iterated logarithm we have that $g(x)$ satisfies conditions \eqref{eqn: g bound} and \eqref{eqn: h conditions} for all $m>\frac{1}{2}$ and combining with the bound \eqref{eqn: singular eigenvalue spread} gives us for every $n$ sufficiently large
\begin{align}\label{eqn: brownian motion eval bound}
    |\lambda_n-r_n|\lesssim n^\frac{1}{m+2}\sqrt{\log\log n}
\end{align}
Therefore, the condition \eqref{eqn: strong sum conditions} holds whenever $\dfrac{n^{\frac{m+3}{m+2}}}{n^\frac{2m}{m+2}}\to 0$. Thus, we conclude that the perturbed operator $\hat{H}$ has a number rigid eigenvalue process whenever $m>3$. Note that this is an example where we need the full power of Proposition \ref{prop: strong sum condition result} and cannot make use of Corollary \ref{prop: weak sum condition result}\\

Consider now the case when $g(x)=S(x)$ is a $\alpha$-stable centered, symmetric, Levy process for $0<\alpha<2$. By Chover's Law of Iterated Logarithm for sums of stable random variables \cite{Chover} we have for all $\epsilon>0$
\[|S_n|\leq n^{1/\alpha}\log(n)^{1/\alpha+\epsilon}\]
for all $n\in\N$ sufficiently large.

Now in order to interpolate between integer times let $n\leq t\leq n+1$, and note that 
\begin{align}\label{eqn: obvious}
|S_t|\leq |S_n|+\displaystyle\sup_{0\leq u\leq 1}|S_{n+u}-S_n|
\end{align}
By \cite[Proposition 4.1]{MaximalInequality} (see also \cite{MaximalInequalityHistorical} for an earlier version of this inequality) we have
\[\P[\sup_{0\leq u\leq 1}|S_{n+u}-S_n|>a]\leq Ca^{-\alpha}\] for some constant $C$ that does not depend on $a$. 

Using this inequality yields

\[\P[\sup_{0\leq u\leq 1}|S_{n+u}-S_n|>n^{1/\alpha}\log(n)^{1/\alpha+\epsilon}]\leq C\frac{1}{n\log(n)^{1+\alpha\epsilon}}\]

Note that this is summable and so by Borel-Cantelli we are able to conclude that 
\[\sup_{0\leq u\leq 1}|S_{n+u}-S_n|=\mathrm{O}(n^{1/\alpha}\log(n)^{1/\alpha+\epsilon}) \quad \text{a.s.}\]

which when combined with \eqref{eqn: obvious} gives us that 
\[\sup_{0\leq t\leq x} |S_t|=\mathrm{O}(x^{1/\alpha}\log(x)^{1/\alpha+\epsilon})\]

 Now we can use \eqref{eqn: singular eigenvalue spread} to get that for $m>\dfrac{1}{\alpha}$
\[|\lambda_n-r_n| = \mathrm{O}(n^{\frac{2}{\alpha(m+2)}}\log(n)^{1/\alpha+\epsilon})  \]
Then the condition \eqref{eqn: strong sum conditions} holds and we get that $\hat{H}$ is number rigid whenever $m>2+2/\alpha$ (note how this aligns with the case where the perturbation is a Brownian Motion). 
\end{proof}

Like our earlier examples, this setting is of great theoretical interest due to the heavy tails of the stable process which render the standard exponential semigroup approach (see Remark \ref{rmk: exponential semigroup}) impossible. We are thus able to see the main advantage of the pathwise statistic \eqref{eqn: strong sum conditions} which gives us access to relatively simple and convenient tools from perturbation theory to prove the number rigidity of $\hat{H}$.

\section{Proof of New Results}

In order to prove Proposition \ref{prop: strong sum condition result} we will explicitly construct a measurable function that gives us the point count of a bounded Borel set given the outside configuration. The design of this function is partially inspired by some previous rigidity results that used similar functions in the role of the test functions for the variance of linear statistics (see e.g. \cite{BufetovAiry}).

\subsection{Proof of Proposition \ref{prop: strong sum condition result}:} \label{sec: proof of strong result}
First we note that it is easy to see from condition \eqref{eqn: strong sum conditions} and the fact that $r_n\to\infty$ that $r_n+g_n\to\infty$ as well.
Now for $R>0$ we define the $\frac{1}{R}$-Lipschitz function
\begin{align}\label{eqn: cutoff function}
    \phi_R(x):=\frac{(R-x)_+}{R}
\end{align}
where \[(\cdot)_+ = \max(\cdot,0)\]
and let 
\begin{align}\label{eqn: centering}
    a_R:=\sum_{i\geq1}\phi_R(r_i)
\end{align}

Now given a particular realization of the point process $\Pi(\omega)$, which corresponds to a particular realization of the perturbation process $\{g_n(\omega)\}$, define the functional 
\[Y_R(\Pi(\omega)) := \int_\R\phi_R(x)d\Pi(x)-a_R := \sum_{i\geq1}\phi_R(r_i+g_i(\omega))-\phi_R(r_i)\]
and observe that 
\begin{align*}
    |Y_R|\leq\frac{1}{R}\sum_{\{i:\min(r_i+g_i,r_i)\leq R\}}|g_i|
\end{align*}
Now taking the sequence $R=r_n$, with $n$ large enough so that $r_n>0$, and writing $m_n:=\max\{i:\min(r_i+g_i,r_i)\leq r_n\}$ yields
\begin{align}\label{eqn: main Y inequality}
        |Y_{r_n}|\leq\frac{1}{r_n}\sum_{i=1}^{m_n}|g_i|
\end{align}
We also have by definition of $m_n$,
\[r_n\geq\min(r_{m_n}+g_{m_n},r_{m_n})\geq r_{m_n}-|g_{m_n}|\geq r_{m_n}-\sum_{i=1}^{m_n}|g_i|\]
Consequently, we have 
\begin{align}\label{eqn: preflip}
    \frac{r_n}{r_{m_n}}\geq 1-\frac{\sum_{i=1}^{m_n}|g_i|}{r_{m_n}}
\end{align}
Taking the reciprocal of \eqref{eqn: preflip} and the limit as $n\to\infty$ with the assumption \eqref{eqn: strong sum conditions} yields
\begin{align}\label{eqn: ratio bound}
    \frac{r_{m_n}}{r_n}\leq \frac{1}{1-\frac{\sum_{i=1}^{m_n}|g_i|}{r_{m_n}}}\to1
\end{align}

Now returning to the inequality \eqref{eqn: main Y inequality}, we can rewrite it as
\begin{align}\label{eqn: Y inequality ratio rewrite}
    |Y_{r_n}|\leq \frac{\sum_{i=1}^{m_n}|g_i|}{r_{m_n}}\frac{r_{m_n}}{r_n}
\end{align}
Using \eqref{eqn: ratio bound}, the assumption \eqref{eqn: strong sum conditions}, and the positivity of $|Y_{r_n}|$ allows us to conclude
\begin{align}\label{eqn: functional goes to 0}
    |Y_{r_n}(\Pi(\omega))|\to0
\end{align}
In addition, for $B\subset\R$ bounded, $\phi_{r_n}
\to1$ uniformly on $B$ as $r_n\to\infty$. Therefore, defining $Y_R$ on restrictions of the configuration $A\subset\R$ Borel by
\[Y_R(A):=\int_{A}\phi_R(x)d\Pi(x)-a_R\] and using the definition of $Y_R$ gives us
\[Y_R(B^c)=Y_R(\Pi(\omega))-\int_B\phi_R(x)d\Pi(\omega)\]
and so by \eqref{eqn: functional goes to 0}
\[Y_{r_n}(B^c)\to-\Pi(\omega)(B) \text{ a.s. as } r_n\to\infty\]
Since it is simply a convergent sum of Lipschitz functions, $\phi_R$ is measurable and therefore $Y_R$ is also measurable with respect to the outside configuration for every choice of bounded $B\subset\R$, hence $\Pi$ is number rigid completing the proof.

\bibliographystyle{alpha}
\bibliography{Bibliography}

 \end{document}